\documentclass[reqno, oneside]{amsart}
\usepackage{style}
\begin{document}
\title[Wiener-Hopf operators]{Wiener-Hopf operators admit \\ triangular factorization}
\author{R.~V.~Bessonov}

\address{St.Petersburg State University ({\normalfont Universitetskaya nab. 7/9, 199034 St.\,Petersburg, Russia}) and St.Petersburg Department of Steklov Mathematical Institute of Russian Academy of Science ({\normalfont Fon\-tan\-ka 27, 191023, St.Petersburg, Russia})}
\email{bessonov@pdmi.ras.ru}

\thanks{The author is supported by RFBR grant mol\_a\_dk 16-31-60053}

\begin{abstract}
We prove that every positive bounded invertible Wiener-Hopf operator admits triangular factorization. This answers the question posed by L.~Sakhnovich in 1994.   
\end{abstract}

\maketitle

\section{Introduction}
A linear bounded operator $T$ on a Hilbert space $H$ is said to admit a left triangular factorization with respect to a chain $\Lit$ of subspaces in $H$ if there are bounded invertible operators $T_1$, $T_2$ from the nest algebra $\A(\Lit)$ generated by $\Lit$ such that 
\begin{equation}\label{eq8}
T = T_1^* T_2.
\end{equation}
A family $\Lit$ of subspaces in $H$ forms a chain if either $E \subset F$ or $F \subset E$ for any pair of subspaces $E, F \in \Lit$. A bounded operator $A$ on $H$ is called upper-triangular with respect to $\Lit$ if $A E \subset E$ for every $E \in \Lit$. The nest algebra $\A(\Lit)$ consists of all bounded operators upper-triangular with respect to $\Lit$. An operator $T \ge 0$ admitting triangular factorization \eqref{eq8} can always be factorized so that $T_1 = T_2$, see \cite{Lar85}.  

\medskip

The general theory of triangular factorization was developed in book  \cite{GK70} by I.~C.~Goh\-berg and M.~G.~Krein, see also   \cite{GGK93} for a modern exposition and \cite{Sakh94} for a summary of known facts. The famous result by D.~R.~Larson \cite{Lar85} says that every positive bounded invertible operator admits a triangular factorization with respect to a given chain $\Lit$ in $H$ if and only if the chain $\Lit$ is countable. In particular, there exists a bounded invertible operator $T \ge 0$ on $L^2(\R_+)$, $\R_+ = [0, +\infty)$, that does not admit triangular factorization \eqref{eq8} with respect to the continuous chain of subspaces $L^2[0,r]$, $r>0$. As to the author knowledge, no concrete example of such an operator is known.

\medskip

Let $\psi$ be a tempered distribution on the real line, $\R$. The Wiener-Hopf operator~$W_{\psi}$ on $L^2(\R_+)$ is densely defined by
$$
W_{\psi}: f \mapsto \int_{\R_+} \psi(t-s)f(s)\,ds, \qquad t \ge 0,
$$
on smooth functions $f$ with compact support in $(0, +\infty)$, where the integral is understood in the distributional sense. The operator $W_{\psi}$ is positive, bounded, and invertible on $L^2(\R_+)$ if and only if $\psi$ is the Fourier image of a function $w > 0$ on~$\R$ such that $w$, $w^{-1}$ are uniformly bounded on $\R$ except of a set of zero Lebesgue measure. 

\medskip

In 1994, L.~Sakhnovich \cite{Sakh94} asked if every positive bounded invertible operator $W_{\psi}$ on $L^2(\R_+)$ admits triangular factorization with respect to the continuous chain of subspaces $L^2[0,r]$, $r>0$. Below we present three equivalent variants of the affirmative answer to this question. 
\begin{Thm}\label{t1}
Every positive bounded invertible operator $W_{\psi}$ on $L^2(\R_+)$ admits triangular factorization $W_{\psi} = A^* A$, where $A$ is a bounded invertible operator on $L^2(\R_+)$ such that $AL^2[0,r] = L^2[0,r]$.
\end{Thm}
Let $\pw_{[0,r]}$ denote the Paley-Wiener space of entire functions $f$, $f \in L^2(\R)$, such that the Fourier spectrum of $f$ is contained in $[0,r]$. For a weight $w>0$ on $\R$ such that $w$, $w^{-1}$ are uniformly bounded on $\R$, the space $\pw_{[0,r]}$ can be identified with a subspace of the weighted space $L^2(w)$.
\begin{Thm}\label{t2}
Let $w$ be a measurable function on $\R$ such that $c_1 \le w(x) \le c_2$ for some positive constants $c_1$, $c_2$ and almost all $x \in \R$. Then there exists an isometric operator $\F_w: L^2(\R_+) \to L^2(w)$ such that $\F_w L^2[0,r] = \pw_{[0,r]}$ for every $r>0$. 
\end{Thm}
In the next theorem, a Hamiltonian $\Hh$ on $\R_+$ is a measurable matrix-valued mapping taking $t \in \R_+$ into a positive semi-definite $2\times 2$ matrix $\Hh(t)$ with real entries such that the function $\trace \Hh$ belongs to $L^1[0, r]$ for every $r>0$ and does not vanish on a set of positive Lebesgue measure. Each Hamiltonian $\Hh$ generates a self-adjoint differential operator $\Di_\Hh$. The correspondence between Hamiltonians $\Hh$ and spectral measures of the operators $\Di_{\Hh}$ they generate is the main issue of Krein--de Branges spectral theory of canonical Hamiltonian systems \cite{Romanov}.   
\begin{Thm}\label{t3}
Let $w$ be a measurable function on $\R$ such that $c_1 \le w(x) \le c_2$ for some positive constants $c_1$, $c_2$ and almost all $x \in \R$. Then there exists the unique Hamiltonian $\Hh$ on $\R_+$ with $\det \Hh = 1$ almost everywhere on $\R_+$ such that $w\,dx$ is the spectral measure for $\Hh$.    
\end{Thm}
In 2012, L.~Sakhnovich \cite{Sakh12} constructed an example of a non-factorable positive bounded invertible Wiener-Hopf operator. An error in his argument was found and discussed by the author in \cite{B17}, where the factorization problem for Wiener-Hopf operators with real symbols $\psi$ was considered. It appears that the method of \cite{B17} can not give Theorem~\ref{t1} in full generality because it heavily uses the diagonal structure of Hamiltonians arising in the corresponding spectral problem.   

\medskip

The initial idea of M.~G.~Krein was to use triangular factorization methods in the spectral theory of differential equations such as string equation, Dirac system, canonical Hamiltonian systems, etc. Indeed if the orthogonal spectral measure of the corresponding differential operator is nice enough, then triangular factorization methods give precise information about the coefficients of the underlying differential equation. A large number of results of this kind is collected in \cite{Den06}, see also \cite{HM04}. On the other hand, in the setting of Theorem \ref{t1} the distribution $\psi$ is too wild for the usage of standard factorization methods. This is why we first prove Theorem~\ref{t3} and then derive from it Theorems \ref{t1}, \ref{t2}. Our approach is based on a technique developed for the proof of Szeg\H{o}-type theorem for canonical Hamiltonian systems \cite{BD2017}, \cite{B18}.          

\section{Proofs}\label{s2}
As it was indicated above, we will use results related to Szeg\H{o}-type theorem for canonical Hamiltonian systems. The general theory of canonical Hamiltonian systems is discussed in \cite{Romanov}, \cite{Winkler95}, \cite{HSW}. The reader interested in a short introduction to the subject could find all information we will need in Section 2 of \cite{B18}.  

\medskip

Next several lemmas are elementary. They will imply an important bound for the Muckenhoupt $A_2$--characteristic of elements of Hamiltonians $\Hh$ whose spectral measures are ``small perturbations'' of the Lebesgue measure on $\R$. We will then approximate the weight $w$ in Theorem \ref{t3} by such measures and take a limit using a compactness argument.      

\medskip

Consider the set $L^1(\R_+) + L^2(\R_+)$ of functions $f$ on $\R_+$ that can be represented in the form $f = f_1 + f_2$ with $f \in L^1(\R_+)$, $f_2 \in L^2(\R_+)$. As usual, we identify functions coinciding almost everywhere on $\R_+$. The set $L^1(\R_+) + L^2(\R_+)$ is the linear normed space with respect to the norm
$$
\|f\|_{1,2} = \inf\bigl\{\|f_1\|_{L^1(\R_+)} + \|f_2\|_{L^2(\R_+)}:\;\; f = f_1 + f_2\bigr\}.
$$
\begin{Lem}\label{l1}
For $f \in L^1(\R_+) + L^2(\R_+)$, one can find functions $f_1$, $f_2$ such that $f = f_1 + f_2$, $|f_{1,2}| \le |f|$ on $\R_+$, and $\|f_1\|_{L^1(\R_+)} + \|f_2\|_{L^2(\R_+)} \le 4\|f\|_{1,2}$. 
\end{Lem}
\beginpf At first, assume that $f$ is positive. Let $f = f_1 + f_2$ for $f_{1,2}$ such that $\|f_1\|_{L^1(\R_+)} + \|f_2\|_{L^2(\R_+)} \le 2\|f\|_{1,2}$. Consider a decomposition of $f_{1} = f_{11} - f_{12}$ such that $f_{11} \ge 0$, $f_{12} \ge 0$, $f_{11}f_{12} = 0$ almost everywhere on $\R_+$. Let $f_{2} = f_{21} - f_{22}$ be a similar decomposition for $f_{2}$. It is straightforward to check that functions 
$\tilde f_{1} = f_{11} - f_{22}$, $\tilde f_2 = f_{21} - f_{12}$ are nonnegative and satisfy $f = \tilde f_1 + \tilde f_2$. In particular, we have $0 \le \tilde f_{1,2} \le f$ almost everywhere on $\R_+$ and $\|\tilde f_{1}\|_{L^1(\R_+)} \le \|f_{1}\|_{L^1(\R_+)}$, $\|\tilde f_{2}\|_{L^2(\R_+)} \le \|f_{2}\|_{L^2(\R_+)}$. This gives the required decomposition for the positive function $f \in L^1(\R_+) + L^2(\R_+)$. For an arbitrary function $f \in L^1(\R_+) + L^2(\R_+)$, let us represent $f$ in the form $f = f^+ - f^-$, where $f^\pm \ge 0$, $f^+f^- = 0$ almost everywhere on $\R_+$. Then $\|f^{\pm}\|_{1,2} \le \|f\|_{1,2}$ and one can take $f_1 = f_1^+ - f_1^-$, $f_2 = f_2^+ - f_2^-$, where $f^+ = f_1^{+} + f_2^{+}$, $f^- = f_1^{-} + f_2^{-}$ are the compositions of positive functions $f^{\pm} \in L^1(\R_+) + L^2(\R_+)$ constructed in the first part of the proof. Now the lemma follows from the triangle inequality. \qed 

\medskip

\noindent We say that a real-valued function $f$ is locally absolutely continuous on $\R_+$ if 
$$
f(t) = c + \int_{0}^{t} f_1(s)\,ds, \qquad t \ge 0,
$$ 
for a constant $c \in \R$ and a function $f_1$ on $\R_+$ such that $f_1 \in L^1[0,r]$ for every $r>0$. Following \cite{BD2017}, define the class $A_2(\R_+, \ell^1)$ to be the set of functions $f \ge 0$ on the half-axis $\R_+$ such that the characteristic 
$$
[f]_{2,\ell^1} = \sum_{n=0}^\infty \left(\int_{n}^{n+2}f(t)\,dt  \int_{n}^{n+2}\frac{1}{f(t)}\,dt - 4\right)
$$ 
is finite. Note that $[f]_{2,\ell^1} \ge 0$ by H\"older's inequality. Since the intervals $[n, n+2]$ and $[n+1, n+3]$ overlap for every $n \ge 0$, we have  $[f]_{2,\ell^1} = 0$ if and only if $f$ is a nonzero constant. This plays an important role for considerations in \cite{BD2017}. Below we will use a different simple feature of the class $A_2(\R_+,\ell^1)$: if $c = \sup_{y>0}[D_{y}f]_{2,\ell^1}$ is finite for $D_y f: t \mapsto f(t/y)$, then $f$ belongs the classical Muckenhoupt class 
$A_2(\R_+)$ and its Muckenhoupt characteristic $[f]_{2}$ is controllable by $c$. Recall that $A_2(\R_+)$ consists of measurable functions $f \ge 0$ on $\R_+$ such that 
\[
[f]_{2} = \sup_{I \subset \R_+}\frac{1}{|I|}\int_{I}f(t)\,dt \cdot \frac{1}{|I|}\int_{I}\frac{1}{f(t)}\,dt < +\infty,
\]
where the supremum is taken over all intervals $I \subset \mathbb{R_+}$.

\medskip

\begin{Lem}\label{l2}
Let $g$ be a positive locally absolutely continuous function such that $g'/g \in L^1(\R_+) + L^2(\R_+)$. Suppose that
$gh + (gh)^{-1} - 2 \in L^1(\R_+)$ for a positive function $h$ on $\R_+$. Then $h \in A_2(\R_+, \ell^1)$ and, moreover, 
$$
[h]_{2,\ell^1} \le c\|gh + (gh)^{-1} - 2\|_{L^1(\R_+)}^{2} + c,
$$ for a constant $c$ depending only on $\|g'/g\|_{1,2}$. 
\end{Lem} 
\beginpf The quantities $[h]_{2,\ell^1}$, $\|g'/g\|_{1,2}$ are invariant with respect to multiplication of $h$, $g$ by a non-zero constant, hence we may assume that $g(0) = 1$. Put $\phi = g'/g$. By construction, we have $\phi \in L^1(\R_+) + L^2(\R_+)$ and 
$g(t) = e^{\int_{0}^{t}\phi(s)\,ds}$, $t \ge 0$.
For $t \in [n, n+1)$, define $h_n: t \mapsto g(n) h(t)$. Let us first show that 
\begin{equation}\label{eq2}
\sum_{n \ge 0} \left(\int_{n}^{n+1} h_n + \int_{n}^{n+1}\frac{1}{h_n} - 2\right) \le c\|gh + (gh)^{-1} - 2\|_{L^1(\R_+)} + c,
\end{equation}
for a constant $c$ depending only on $\|\phi\|_{1,2}$. Denote by $N_1$ the set of indexes $n$ such that $\int_{n}^{n+1}|\phi(t)|\,dt < 1/4$. For $n \in N_1$ and $t \in [n, n+1)$, we have 
\begin{align}
\sqrt{g(t)h(t)} - \frac{1}{\sqrt{g(t)h(t)}} 
&= e^{\frac{1}{2}\int_{n}^{t}\phi(s)\,ds}\sqrt{h_n(t)} - e^{-\frac{1}{2}\int_{n}^{t}\phi(s)\,ds}\frac{1}{\sqrt{h_n(t)}}\notag \\
&= \sqrt{h_n(t)} - \frac{1}{\sqrt{h_n(t)}} + E_n(t),\label{eq9}
\end{align}
where $|E_n(t)| \le 2\int_{n}^{n+1}|\phi(s)|\,ds \cdot \bigl(\sqrt{h_n(t)} + 1/\sqrt{h_n(t)}\bigr)$. Observe that
\begin{align}
\int_{n}^{n+1}|E_n(t)|^2\,dt 
&\le 4\left(\int_{n}^{n+1}|\phi(t)|\,dt \right)^2 \left(\int_{n}^{n+1}\Bigl(\sqrt{h_n(t)} - 1/\sqrt{h_n(t)}\Bigr)^2\,dt + 4\right) \notag\\
&\le \frac{1}{4}\int_{n}^{n+1}\Bigl(\sqrt{h_n(t)} - 1/\sqrt{h_n(t)}\Bigr)^2 + 16c_n^2, \label{eq10}
\end{align}
where $c_n = \int_{n}^{n+1}|\phi(s)|\,ds$. From \eqref{eq9} we get
$$
\int_{n}^{n+1}\Bigl(\sqrt{h_n} - 1/\sqrt{h_n}\Bigr)^2\,dt 
\le 2 \int_{n}^{n+1}\Bigl(\sqrt{gh} - 1/\sqrt{gh}\Bigr)^2\,dt + 2\int_{n}^{n+1}|E_n|^2\,dt,
$$
so \eqref{eq10} implies
\begin{equation}\label{eq11}
\frac{1}{2}\int_{n}^{n+1}\Bigl(\sqrt{h_n} - 1/\sqrt{h_n}\Bigr)^2\,dt \le 
2 \int_{n}^{n+1}\Bigl(\sqrt{gh} - 1/\sqrt{gh}\Bigr)^2\,dt + 32c_n^2.
\end{equation}
Let $\phi_1 \in L^1(\R_+)$, $\phi_2 \in L^2(\R_+)$ be the functions given by Lemma \ref{l1} for $\phi$. Using the fact that $\int_{n}^{n+1}|\phi_1(s)|\,ds \le c_n \le 1$ for $n \in N_1$, we obtain 
\begin{align*}
\sum_{n \in N_1} c_n^2 
&\le 2 \sum_{n \in N_1}\left(\int_{n}^{n+1}|\phi_1(s)|\,ds\right)^2 + 2\sum_{n \in N_1}\int_{n}^{n+1}|\phi_2(s)|^2\,ds,  \\
&\le 2\|\phi_{1}\|_{L^1(\R_+)} + 2\|\phi_{2}\|_{L^2(\R_+)}^{2}. 
\end{align*} 
From \eqref{eq11} we now see that
$$
\sum_{n \in N_1}\int_{n}^{n+1}(h_n + h_n^{-1} -2)\,ds \le 4\|gh + (gh)^{-1} -2\|_{L^1(\R_+)} + 64(\|\phi\|_{1,2} + \|\phi\|_{1,2}^{2}). 
$$
Next, consider the set $N_2$ of indexes $n$ such that $\int_{n}^{n+1}|\phi(t)|\,dt \ge 1/4$. There are at most $8\|\phi_1\|_{L^1(\R_+)}$ integers $n \ge 0$ such that $\int_{n}^{n+1}|\phi_1(t)|\,dt \ge 1/8$ and at most $64\|\phi_2\|_{L^2(\R_+)}^{2}$ integers $n \ge 0$ such that $\int_{n}^{n+1}|\phi_2(t)|\,dt \ge 1/8$.
It follows that the number of elements in $N_2$ does not exceed $64(\|\phi\|_{1,2} + \|\phi\|_{1,2}^{2})$. Using the trivial bound 
\begin{align}
\int_{n}^{n+1}\left(h_n + h_{n}^{-1}\right)\,dt 
&\le e^{\int_{n}^{n+1}|\phi(s)|\,ds} \int_{n}^{n+1}\left(gh + (gh)^{-1}\right)\,dt \notag \\
&\le e^{\|\phi\|_{1,2}}\Bigl(\|gh + (gh)^{-1} -2\|_{L^1(\R_+)} + 2\Bigr), \label{eq3}
\end{align}
for $n \in N_2$, we conclude that \eqref{eq2} holds. Next, for $n \ge 0$ put
$$
v_n = \int_{n}^{n+1}h, \quad w_n = \int_{n}^{n+1}h^{-1},
\quad
\tilde v_n = \int_{n}^{n+1}h_n, \quad \tilde w_n = \int_{n}^{n+1}h_{n}^{-1}.
$$
We have 
\begin{equation}\label{eq41}
\sum_{n \ge 0} \left(2\sqrt{v_n w_n} - 2\right) = \sum_{n \ge 0} \left(2\sqrt{\tilde v_n  \tilde w_n} - 2\right) 
\le \sum_{n \ge 0} \left(\tilde v_n + \tilde w_n - 2\right).
\end{equation}
We also have $v_n w_n = \tilde v_n \tilde w_n \le e^{2\|\phi\|_{1,2}}\Bigl(\|gh + (gh)^{-1} -2\|_{L^1(\R_+)} + 2\Bigr)^2$
by \eqref{eq3}, hence the inequalities \eqref{eq2}, \eqref{eq41}, and 
$$
\sum_{n \ge 0} (v_n w_n - 1) \le (\max_{n \ge 0}\sqrt{v_n w_n} + 1)\sum_{n \ge 0} (\sqrt{v_n w_n} - 1),
$$
yield
$$
\sum_{n \ge 0}\left(\int_{n}^{n+1}h(t)\,dt\int_{n}^{n+1}\frac{dt}{h(t)} - 1\right) \le c\|gh + (gh)^{-1} - 2\|_{L^1(\R_+)}^{2} + c,
$$
for a constant $c$ depending only on $\|g'/g\|_{1,2}$. Similar estimates hold for the pairs of functions $t \mapsto g(2t)$, $t \mapsto h(2t)$, and $t \mapsto g(2t - 1)$, $t \mapsto h(2t - 1)$, that satisfy all assumptions of the lemma. These estimates imply the desired bound for $[h]_{2, \ell^1}$. \qed

\medskip

Below we follow notations and definitions from papers \cite{BD2017}, \cite{B18}. In particular, the spectral measure $\mu$ of a Hamiltonian $\Hh$ is the measure in the Herglotz representation of its Weyl function
\begin{equation}\label{eq14}
m(z) = \lim_{t \to +\infty}\frac{\Phi^-(t,z)}{\Theta^-(t,z)} 
= \frac{1}{\pi}\int_{\R}\left(\frac{1}{x-z} - \frac{x}{1+x^2}\right)d\mu(x) + bz + a,
\end{equation}  
in the upper half-plane $\C^+ = \{z \in \C: \Im z > 0\}$. The functions $\Phi^-$, $\Theta^-$ above are the entries of the solution $M = \mms$ of Cauchy problem 
\begin{equation}\label{cs}
\begin{cases}
J \tfrac{\partial}{\partial t} M(t,z) = z \Hh(t) M(t,z), \\ M(0, z) = \idm,
\end{cases} 
\quad t \in \R_+, \quad z \in \C,
\end{equation}
where $J = \jm$. The singular part of a measure $\mu$ on $\R$ will be denoted by $\mus$. The Szeg\H{o} class $\sz$ consists of measures $\mu = w\,dx + \mus$ on $\R$ such that $(1+x^2)^{-1} \in L^1(\mu)$ and $(\log w(x))/(1+x^2) \in L^1(\R)$. For such measures $\mu$ we will use notation 
$$
\K(\mu, z) = \log\frac{1}{\pi}\int_{\R}\frac{\Im z}{|x - z|^2}\,d\mu(x) - \frac{1}{\pi}\int_{\R}\log w(x)\frac{\Im z}{|x - z|^2}\,dx,
\qquad z \in \C^+.
$$   
We also put $\K(\mu) = \K(\mu, i)$. By Jensen inequality, $\K(\mu, z) \ge 0$ for every $z \in \C^+$. For Hamiltonians $\Hh$ such that $\det \Hh \neq 0$ almost everywhere on $\R_+$, we have 
$b = 0$ in \eqref{eq8} by Lemma 2.3 in \cite{BD2017}. In particular, for the spectral measure $\mu$ of such a Hamiltonian, $\K(\mu)$ coincides with the quantity $\K_\Hh(0)$ defined in Section 2.2 of \cite{B18}. This observation and Lemma 3 in \cite{B18} yield the following result.
\begin{Lem}\label{l9}
Let $\Hh$ be a Hamiltonian on $\R_+$ such that $\det \Hh \neq 0$ almost everywhere on $\R_+$, and let $\Hh^d = J^* \Hh J$. Assume that the spectral measure $\mu = w(x)\,dx + \mus$ of $\Hh$ is such that $\mu \in \sz$. Then $\mu^d \in \sz$ and $\K(\mu^d) = \K(\mu)$ for the spectral measure $\mu^d$ of the Hamiltonian $\Hh^d$. 
\end{Lem}
Another result we will need is a weak variant of the Szeg\H{o} theorem \cite{BD2017} for canonical Hamiltonian systems.
\begin{Lem}\label{l8}
Let $\Hh = \sth$, $\mu$ be as in Lemma \ref{l9}. Then there exist positive locally absolutely continuous functions $g_1$, $g_2$ on $\R_+$ such that 
\begin{align*}
\left\|g_k h_k + (g_k h_k)^{-1} - 2\right\|_{L^1(\R_+)} &\le c_\mu, \qquad \|g'_k/g_k\|_{1,2} \le c_\mu, 
\end{align*}
for $k = 1,2$, and a constant $c_\mu$ controllable by $\K(\mu)$.
\end{Lem}
\beginpf For $k =1$, the statement is a corollary of the proof Lemma 7 in \cite{B18} (see formula (39) therein). For $k=2$, one needs to consider the dual Hamiltonian $\Hh^d = \sthd$ and use the fact that $\K(\mu) = \K(\mu^d)$ from Lemma \ref{l9}. \qed 

\medskip

In the next lemma we obtain the key estimate for what follows.
\begin{Lem}\label{l3}
Let $\Hh = \sth$ be a Hamiltonian on $\R_+$ such that $\det \Hh = 1$ almost everywhere on $\R_+$. Assume that the spectral measure $\mu = w(x)\,dx + \mus$ of $\Hh$ is such that $\sup_{y>0} \K(\mu, iy) < +\infty$. Then $h_1$, $h_2$ belong to $A_2(\R_+)$ and, moreover, we have $[h_{1}]_{2} \le c_\mu$, $[h_{2}]_{2} \le c_\mu$, for a constant $c_\mu$ depending only on $\sup_{y>0} \K(\mu, iy)$.  
\end{Lem}
\beginpf For $y >0$, consider the Hamiltonian $D_y \Hh: t \mapsto \Hh(t/y)$ on $\R_+$. By construction, $\det D_y\Hh = 1$ almost everywhere  on $\R_+$. If $M(t,z)$ is the solution of Cauchy problem \eqref{cs}, then $t \mapsto M(t/y, zy)$ is the solution of the same Cauchy problem for the Hamiltonian $D_y\Hh$. It follows that the Weyl function $m^y$ of $D_y \Hh$ is given by
\begin{equation}\label{eq12}
m^{y}(z) = \lim_{t \to +\infty}\frac{\Phi_{-}(t/y,yz)}{\Theta_{-}(t/y,yz)} = m(yz), \qquad z \in \C^+, 
\end{equation}
where $m$ is the Weyl function for $\Hh$. Denoting by $\mu^y = w^{y}\,dx + \mus^y$ the spectral measure of $\Hh^y$, from  \eqref{eq12} we see that 
$$
\frac{1}{\pi}\int_{\R}\frac{d\mu^y(x)}{x^2 + 1} = \Im m^y(i) = \Im m(iy) = \frac{1}{\pi}\int_{\R}\frac{y}{x^2 + y^2}\,d\mu(x).
$$ 
Here we used the fact that $b = 0$ in \eqref{eq14} for Hamiltonians $\Hh$ with $\det\Hh \neq 0$ on $\R_+$, see Lemma 2.3 in \cite{BD2017}. We also have 
$$
w^y(x) = \lim_{\eps \to +0} \Im m^y(x+i\eps)= \lim_{\eps \to +0} \Im m(yx+iy\eps) = w(yx),
$$ 
for almost all $x \in \R$.  It follows that $\K(\mu^y) = \K(\mu, iy)$ and $\sup_{y>0} \K(\mu^y) < +\infty$. By Lemma \ref{l8}, there exists a positive locally absolutely continuous function $g$ on $\R_+$ and a constant $\tilde c_\mu$ depending only on $\sup_{y>0}\K(\mu, iy)$ such that $g'/g \in L^1(\R_+) + L^2(\R_+)$ with $\|g'/g\|_{1,2} \le \tilde c_\mu$, and 
$$
\left\|g D_y h_1 + (g D_y h_1)^{-1} - 2\right\|_{L^1(\R_+)} \le \tilde c_\mu, \qquad D_y h_1: t \mapsto h_1(t/y). 
$$
From Lemma \ref{l2} we see that $0 \le [D_y h_1]_{2, \ell^1} \le c_\mu - 4$ for another constant $c_\mu$ depending only on $\sup_{y>0}\K(\mu, iy)$, where we subtract $4$ for the future convenience. In particular, we have
$$
\int_{n}^{n+2}D_y h_1(t)\,dt \cdot \int_{n}^{n+2}\frac{dt}{D_y h_1(t)} \le c_\mu,
$$ 
for every integer $n \ge 0$ and every $y>0$. This can be rewritten in the form 
\begin{equation}\label{eq51}
\sup_{I}\frac{1}{|I|}\int_{I}h_1(t)\,dt \cdot \frac{1}{|I|}\int_{I}\frac{dt}{h_1(t)} \le c_\mu/4,
\end{equation} 
where the supremum is taken over all intervals $I$ of the form $I = [n/y, (n+2)/y]$. It is easy to see that for every interval $J \subset \R_+$ one can find interval $I \supset J$ of this form such that $|I| \le 2|J|$. Hence the supremum in \eqref{eq51} over all intervals of $\R_+$ does not exceed $c_\mu$, that is, $[h_{1}]_2 \le c_\mu$. The same consideration applies to the Hamiltonian $\Hh^d = \sthd$, yielding $[h_{2}]_2 \le c_\mu$. \qed 

\medskip

Given a Hamiltonian $\Hh$ on $\R_+$ with the spectral measure $\mu$, denote by $(\pw_{[0,r]}, \mu)$ the Hilbert space of entire functions
$$
(\pw_{[0,r]}, \mu) = \left\{f: \; f = e^{irz/2}\tilde f \; \mbox{ for } \tilde f \in \B_r\right\},
$$
with the inner product inherited from $L^2(\mu)$. Here $\B_r$ is the de Branges space generated by the restriction of $\Hh$ to the interval $[0,r]$, see Section 2.3 in \cite{B18} for precise definition. In the case where $\mu = w(x)\,dx$ for a measurable function $w \ge 0$ such that $w$, $w^{-1}$ are uniformly bounded on $\R$, the space $(\pw_{[0,r]}, \mu)$ coincides as a set with the usual Paley-Wiener space $\pw_{[0,r]}$ defined in the Introduction. Consideration of this particular case is sufficient for the proof the main results of the paper, but we will treat the general situation in the next lemma.
  
\medskip

\begin{Lem}\label{l4}
Let $\Hh$ be a Hamiltonian on $\R_+$ such that $\det \Hh = 1$ 
almost everywhere on $\R_+$, and let $\mu$ be its spectral measure. Then there exist entire functions $\{P_{r}\}_{r>0}$ such that for every $r>0$ the mapping 
\begin{equation}\label{eq7}
\F_{\mu}: f \mapsto \frac{1}{\sqrt{2\pi}}\int_{0}^{r}f(t)P_{t}(z)\,dt, \qquad z \in \C,
\end{equation} 
is the unitary operator from $L^2[0,r]$ into $(\pw_{[0,r]}, \mu)$ densely defined on functions in $L^2(\R_+)$ with compact support. 
\end{Lem}
\beginpf Let $\Theta = \ths$ be the first column of the solution $M$ of Cauchy problem~\eqref{cs}. Choose a representative of the mapping $\sqrt{\Hh}$ and define for $r \ge 0$, $z \in \C$,
$$
\Psi = \left(\begin{smallmatrix}\Psi^+ \\ \Psi^-\end{smallmatrix}\right) =\sqrt{\Hh}\Theta, \qquad P_{2r}(z) = e^{irz}(\Psi^+(r,z) - i \Psi^-(r,z)). 
$$ 
We have
\begin{align*}
|P_{2r}(z)|^2 
&= e^{-2r\Im z}\|\Psi(r,z)\|_{\C^2}^{2} - 2e^{-2r\Im z} \Im(\Psi^+(r,z)\ov{\Psi^-(r,z)})\\
&= e^{-2r\Im z}\|\Psi(r,z)\|_{\C^2}^{2}  + i e^{-2r\Im z} \left\langle J \Psi(r,z), \Psi(r,z)\right\rangle_{\C^2}\\
&= e^{-2r\Im z}\left\langle \Hh(r)\Theta(r,z), \Theta(r,z)\right\rangle_{\C^2} + i e^{-2r\Im z} \left\langle J\Theta(r,z), \Theta(r,z)\right\rangle_{\C^2},
\end{align*}
where we used the fact that $A J A = J$ for every real matrix $A = A^*$ with unit determinant, $J = \jm$. The well-known identity
$$
\left\langle J\Theta(r,z), \Theta(r,z)\right\rangle_{\C^2} = 2i \Im z \int_{0}^{r}\left\langle \Hh(t)\Theta(t,z), \Theta(t,z)\right\rangle_{\C^2}\,dt, \qquad r \ge 0,  
$$ 
simply follows from \eqref{cs} by differentiation. We now see that 
$$
|P_{2r}(z)|^2 = \frac{\partial}{\partial r}\left(e^{-2r\Im z} \int_{0}^{r}\left\langle \Hh(t)\Theta(t,z), \Theta(t,z)\right\rangle_{\C^2}\,dt\right), 
$$
for almost every $r \ge 0$. Thus, for all $r\ge 0$, $z \in \C$, we have
\begin{align*}
\int_{0}^{r}|P_{t}(z)|^2\,dt 
&= 2\int_{0}^{r/2}|P_{2t}(z)|^2\,dt \\
&= 2e^{-r\Im z} \int_{0}^{r/2}\left\langle \Hh(t)\Theta(t,z), \Theta(t,z)\right\rangle_{\C^2}\,dt \\
&= 2e^{-r\Im z} \frac{\left\langle J\Theta(r/2,z), \Theta(r/2,z)\right\rangle_{\C^2}}{z - \bar z}.
\end{align*}
In particular, the function $t \mapsto P_t(z)$ is in $L^{2}_{\rm loc}(\R_+)$, and the integral in \eqref{eq7} converges for functions $f \in L^2(\R_+)$ with compact support. Moreover, for all $z, \lambda \in \C$ we have
\begin{align}
\int_{0}^{r}P_{t}(z)\ov{P_{t}(\lambda)}\,dt 
&= 2 e^{ir(z - \bar\lambda)/2} \frac{\left\langle J\Theta(r/2,z), \Theta(r/2,\lambda)\right\rangle_{\C^2}}{z - \bar \lambda},\label{eq5}\\
&= 2\pi e^{ir(z - \bar\lambda)/2} \left(-\frac{1}{2\pi i}\frac{E(z)\ov{E(\lambda)} - E^\sharp(z)\ov{E^\sharp(\lambda)}}{z - \bar \lambda}\right), \notag
\end{align}
where $E(z) = \Theta^+(r/2,z) + i \Theta^-(r/2,z)$, and $E^\sharp(z) = \Theta^+(r/2,z) - i \Theta^-(r/2,z)$. The right hand side of the above identity coincides with $2\pi k_{r,\lambda}(z)$, where $k_{r, \lambda}$ is the reproducing kernel of the space $(\pw_{[0, r]}, \mu)$, see Section 2.3 in \cite{B18}. For $r > 0$ and $z \in \C$, denote by $e_{r, \lambda}$ the function 
$t \mapsto \chi_{[0, r]}(t) \ov{P_{t}(\lambda)}$, where $\chi_{[0,r]}$ is the indicator function of the interval $[0,r]$. From formula \eqref{eq5} we see that 
$\F_{\mu}e_{r,\lambda} = \sqrt{2\pi} k_{r, \lambda}$, and, moreover, 
$$
(\F_{\mu}e_{r,\lambda}, \F_{\mu}e_{r,z})_{L^2(\mu)} = 2\pi(k_{r,\lambda}, k_{r,z})_{L^2(\mu)} = 2\pi k_{r, \lambda}(z) = (e_{r,\lambda}, e_{r,z})_{L^2(\R_+)},
$$
where we used the reproducing kernel property $f(z) = (f, k_{r,z})_{L^2(\mu)}$ for the function $f = k_{r,\lambda}$ in $(\pw_{[0,r]}, \mu)$. In other words, the operator $\F_\mu: L^2(\R_+) \to L^2(\mu)$ is correctly defined and isometric on the linear span of $\{e_{r,\lambda},\; r > 0,  \lambda \in \C\}$. Observe that this linear span is dense~$L^2(\R_+)$. Indeed if a compactly supported function $f \in L^2(\R_+)$ is orthogonal to $e_{r, 0}$ for all $r>0$, then
$$
\int_{0}^{r} f(t)P_t(0)\,dt =0, \qquad r > 0.
$$ 
Hence $f(t)P_t(0) = 0$ almost everywhere on $\R_+$. By construction, 
$$P_{2t}(0) = \Psi^+(t,0) - i \Psi^-(t,0) = a_{11}(t) - ia_{21}(t),$$
where $a_{11}$, $a_{21}$ are the entries of $\sqrt{\Hh}$. Since $\sqrt{\Hh}$ is real and $\det \sqrt{\Hh} = 1$ almost everywhere on $\R_+$, we cannot have $a_{11}(t) - ia_{21}(t) = 0$ on a subset of $\R_+$ of positive Lebesgue measure. Hence $f = 0$ and the set $\{e_{r,\lambda},\; r > 0, \,  z \in \C\}$ is complete in $L^2(\R_+)$. Thus, the operator $\F_{\mu}$ is isometric from $L^2(\R_+)$ to $L^2(\mu)$. Since it sends the elements $e_{r,\lambda} \in L^2[0,r]$ into the complete family $\{k_{r, \lambda}, \lambda \in \C\}$ of reproducing kernels of $(\pw_{[0,r]}, \mu)$, we have $\F_\mu L^2[0,r] = (\pw_{[0,r]}, \mu)$, as claimed. \qed   

\bigskip

\noindent {\bf Proof of Theorem \ref{t3}.} Let $w$ be a function on $\R_+$ such that $0< c_1 \le w(x) \le c_2$ for some constants $c_{1}$, $c_2$ and almost all $x \in \R$. Consider the sequence of functions $w_j$ defined by 
$$
w_j(x) = \begin{cases}w(x), &|x| \le j\\ 1, &|x|>j\end{cases}.
$$
Since the Fourier transform of $1 - w_j$ is infinitely smooth, there exists a smooth Hamiltonian $\Hh_j$ such that 
$\mu_j = w_j(x)\,dx$ is the spectral measure for $\Hh_j$, and $\det\Hh_j = 1$ almost everywhere on $\R_+$. This known fact follows from the classical Gelfand-Levitan theory for Dirac systems, see, e.g., Lemma 3.2 in \cite{B17}. Observe that $\mu_j \in \sz$ for every $j \ge 0$, and the quantities $\sup_{y>0}\K(\mu_j, iy)$ are uniformly bounded in $j \ge 1$. By Lemma \ref{l3}, there exists a constant $c$ such that
$$
[h_{j,1}]_{2} \le c, \qquad [h_{j,2}]_{2} \le c, \qquad \Hh_j = \left(\begin{smallmatrix}h_{j,1} & h_j \\ h_j & h_{j,2}\end{smallmatrix}\right), 
$$
for all $j \ge 1$. We next proceed as in the proof of Theorem 1 in \cite{B17}. Fix $r > 0$. Using the ``reverse H\"older's inequality'' for Muckenhoupt weights \cite{Vas03}, we see that there are constants $c(r)$, $p>1$, depending only on $r$ and $\sup_{j}([h_{j,1}]_{2} + [h_{j,2}]_{2})$ such that $h_{j,k} \in L^p[0, r]$ and $\|h_{j,k}\|_{L^p[0,r]} \le c(r)$ for $k = 1,2$. Since $h_{1,j}h_{2,j}- h_j^2 = \det \Hh_{j} = 1$, we also have
$$
\left(\int_{0}^{r}|h_j(t)|^p \,dt \right)^{1/p} \le \left(\left(\int_{0}^{r}(h_{j,1}(t)h_{j,2}(t))^{\frac{p}{2}}\,dt\right)^{2/p} \!+ r^{2/p}\right)^2 \le 
\left(c(r)^2 + r^{2/p}\right)^{2}. 
$$
Using the diagonalization procedure, one can find subsequences $h_{j_{k}, 1}$, $h_{j_{k}, 2}$, $h_{j_{k}}$ converging weakly in $L^p[0,r]$ for every $r>0$. Let $h_1$, $h_2$, $h$ denote the weak limits of $h_{j_{k}, 1}$, $h_{j_{k}, 2}$, $h_{j_{k}}$, correspondingly. Repeating literally the argument from the proof of Theorem 1 in \cite{B17}, we see that $\Hh = \sth$ is the Hamiltonian on $\R_+$ such that $\det \Hh = 1$ almost everywhere on $\R_+$, and $\mu = w\,dx$ is the spectral measure for~$\Hh$. \qed 

\medskip

\noindent {\bf Proof of Theorem \ref{t2}.} The result is the direct consequence of Theorem \ref{t3} and Lemma \ref{l4}. \qed 

\medskip

\noindent {\bf Proof of Theorem \ref{t1}.} Let $\F : f \mapsto \frac{1}{\sqrt{2\pi}}\int_{\R}f(x)e^{-ixt}\,dx$ denote the usual Fourier transform on $L^2(\R)$. Since the operator $W_{\psi}$ on $L^2(\R_+)$ is positive, bounded and invertible, the same is true for the Toeplitz operator $T = \F W_{\psi}\F^{-1}$ on the Hardy space $H^2 = \F L^2(\R_+)$. Hence the symbol $w$ of the operator $T = T_w$ is such that $0<c_1 \le w(x) \le c_2$ for some constants $c_1$, $c_2$ and almost all $x \in \R$, see Section 4.2.7 in Part B of \cite{Nik02t1}. Denote $\mu = w\,dx$ and consider the isometric mapping $\F_{\mu}: L^2(\R_+) \to L^2(\mu)$ from Theorem \ref{t2}. We have
$\F_{\mu}L^2[0,r] = (\pw_{[0, r]}, \mu)$ for every $r>0$. Let $H^2(\mu) = \F_{\mu} L^2(\R_+)$ denote the weighted Hardy space in $\C^+$ with the inner product inherited from $L^2(\mu)$. Since $w$, $w^{-1}$ are uniformly bounded, the identity embedding $j : H^2 \to H^2(\mu)$ is bounded, invertible, and such that $j \pw_{[0,r]} = (\pw_{[0, r]}, \mu)$. Taking $A = \F\F_{\mu}^{-1}j$, we see that
$$
(A^*A f, f)_{L^2(\R)} = \|\F\F_{\mu}^{-1}j f\|^2_{L^2(\R)} = \|j f\|^2_{L^2(\R)} = \|f\|^2_{L^2(\mu)} = (T_w f,f). 
$$ 
In other words, $T_w = A^*A$ admits the triangular factorization along the chain of Paley-Wiener subspaces $\{\pw_{[0,r]}\}_{r>0}$ of the space $H^2$. Taking the Fourier transform, we conclude that the Wiener-Hopf operator $W_{\psi} = \F^{-1} T_w \F$ on $L^2(\R_+)$ admits the triangular factorization $W_{\psi} = \hat A^* \hat A$, $\hat A = \F^{-1} A\F$, along the chain $\{L^2[0,r]\}_{r>0}$ of subspaces in $L^2(\R_+)$. The theorem follows. \qed

%\vspace{-0.1cm}

%\bibliographystyle{unsrt}
\bibliographystyle{plain} 
\bibliography{bibfile}

\enddocument